\documentclass[aos]{imsart}

\RequirePackage{amsthm,amsmath}
\usepackage{amsfonts}
\usepackage{amssymb}
\usepackage{graphicx}

\RequirePackage{natbib}

\startlocaldefs
\input pollard.sty

\def\Leb{{\mathfrak m}}
\def\supnl{\sup\nolimits}

\def\subheading#1{\medskip\noindent{\normalfont\bfseries #1}}

\theoremstyle{plain}
\newtheorem{theorem}[equation]{Theorem}
\newtheorem{definition}[equation]{Definition}

\theoremstyle{remark}
\newenvironment{remark}{\begin{quote}{\bf Remark.\quad}\rm\small}{\end{quote}}

\def\Theorem{\begin{theorem}}
\def\endTheorem{\end{theorem}}

\def\Remark{\begin{remark}}
\def\endRemark{\end{remark}}

\endlocaldefs

\begin{document}

\begin{frontmatter}
\title{A note on insufficiency and the preservation of Fisher information}
\runtitle{No loss of information}

\begin{aug}
\author{
\fnms{David} 
\snm{Pollard}
\ead[label=e1]{david.pollard@yale.edu}
\ead[label=u1,url]{http://www.stat.yale.edu}}

\runauthor{David Pollard}

\affiliation{Yale University}

\address{Statistics Department\\
Yale University\\
New Haven, Connecticut\\
}

\end{aug}

\begin{abstract}
\citet{KaganShepp05AmStat} presented an elegant example of a mixture model  for which an insufficient statistic  preserves Fisher information.
This note uses the regularity property of differentiability in quadratic mean to provide another explanation for the phenomenon they observed. Some connections with Le~Cam's theory for convergence of experiments are noted. 
\end{abstract}



\end{frontmatter}

\section{Introduction}  \label{intro}
Suppose $\pp=\{P_\th: \th\in\Th\}$ is a
statistical experiment, a set of probability measures
on some measure space~$(\xx,\aa)$ indexed by a subset~$\Th$ of the real line.

The Fisher information function~$\II_\pp(\th)$ can be defined under
various regularity conditions.  If $S$ is a measurable map from  $\xx$
into another measure space $(\yy,\bb)$, each image measure $Q_\th =
S P_\th$ (often called the distribution of~$S$ under~$P_\th$, and sometimes denoted
by~$P_\th S^{-1}$) is a probability measure on~$\bb$. The
statistical experiment $\qq=\{Q_\th:\th\in\Th\}$ is  less
informative, in the sense that an observation~$y\sim Q_\th$ tells us
less about~$\th$ than an observation~$x\sim P_\th$.  In particular,
$\II_\qq(\th)\le \II_\pp(\th)$ for every~$\th$. If $S$ is a sufficient
statistic the last inequality becomes an equality: there is no loss of
Fisher information.

Statistical folklore holds that the converse is also true. For example, \citet[page~158]{LehmannCasella1998}
set as an exercise the task of verifying, ``under suitable regularity conditions'',
results stated by \citet[Section~1]{Basu1964Ancillary}, including the assertion that  there is no loss of Fisher information if and only if the statistic is sufficient. They interpreted Basu's (unstated) regularity conditions to be ``mainly concerned with interchange of integration with differentiation". 
Nevertheless,
\citet{KaganShepp05AmStat} (henceforth K\&S) were able to show, by means of a
simple example, that it is possible to have $\II_\qq(\th)=
\II_\pp(\th)$ for every~$\th$ without~$S$ being sufficient. 
The K\&S counterexample relies on another property---the support of a density changing with the unknown parameter---that is notorious for upsetting classical statistical theory.

My purpose in this note is to make two small additions to the K\&S analysis. 
\begin{enumerate}
\item 
I  reinterpret the phenomenon
identified by K\&S, using the geometry of
differentiabilty in quadratic mean.
\item
Using an asymptotic argument,
I  offer an explanation for why the extent of the failure of sufficiency in the K\&S example is too small to be captured by the Fisher information. More precisely, I  explain why the experiment~$\qq_n$
obtained by~$n$ independent replications of~$\qq$ is asymptotically
equivalent (in Le~Cam's sense) to the corresponding~$\pp_n$.
\end{enumerate}

Most of the necessary theory is already available in the literature
but is not widely known.  The K\&S example provides a good showcase 
for that theory.

\section{The K\&S example} \label{KSeg}
What follows is a slightly simplified version of the K\&S
construction.

Start from a smooth probability density
$$
g(w) = \tfrac12 w^2 e^{-w}\{w>0\}
$$
\wrt/ Lebesgue measure $\Leb$ on the real line. 
The power $w^2$ is chosen so that
$$
\frac{\gdot(w)^2}{g(w)} = g(w) \left(\frac{d\log g(w)}{dw}\right)^2=\frac12 (2-w)^2e^{-w}\{w>0\}
$$
is Lebesgue integrable. The shift family of densities~$\{g(w-\th): \th\in \RR\}$ has constant Fisher information,
\begin{equation}
\II  =  \int_{-\infty}^{\infty} \gdot(w)^2/g(w)\,dw < \infty  .
\label{II.def}
\end{equation}

Let $\nu$ denote the probability measure that puts mass~$1/2$ at each
of~$+1$ and~$-1$. For  each $\th\in\Theta=\RR$ define a probability
measure $P_\th$ on (the Borel sigma-field of) $\xx=\RR\times\{-1,+1\}$
by means of its density
\begin{equation}
f_\th(x) = \{z=+1\}g(y-\th) +\{z= -1\}g(\th-y)
\qt{where $x = (y,z)\in\xx$}
\label{f.def}
\end{equation}
\wrt/~the measure~$\lam:=\Leb\otimes\nu$. That is, the coordinate~$z$
has marginal distribution~$\nu$ and the conditional distribution of
$y$ given~$z$ is that of $\th+zw$ where~$w\sim g$ independently of~$z$.

\medskip
\includegraphics[width=4.7in]{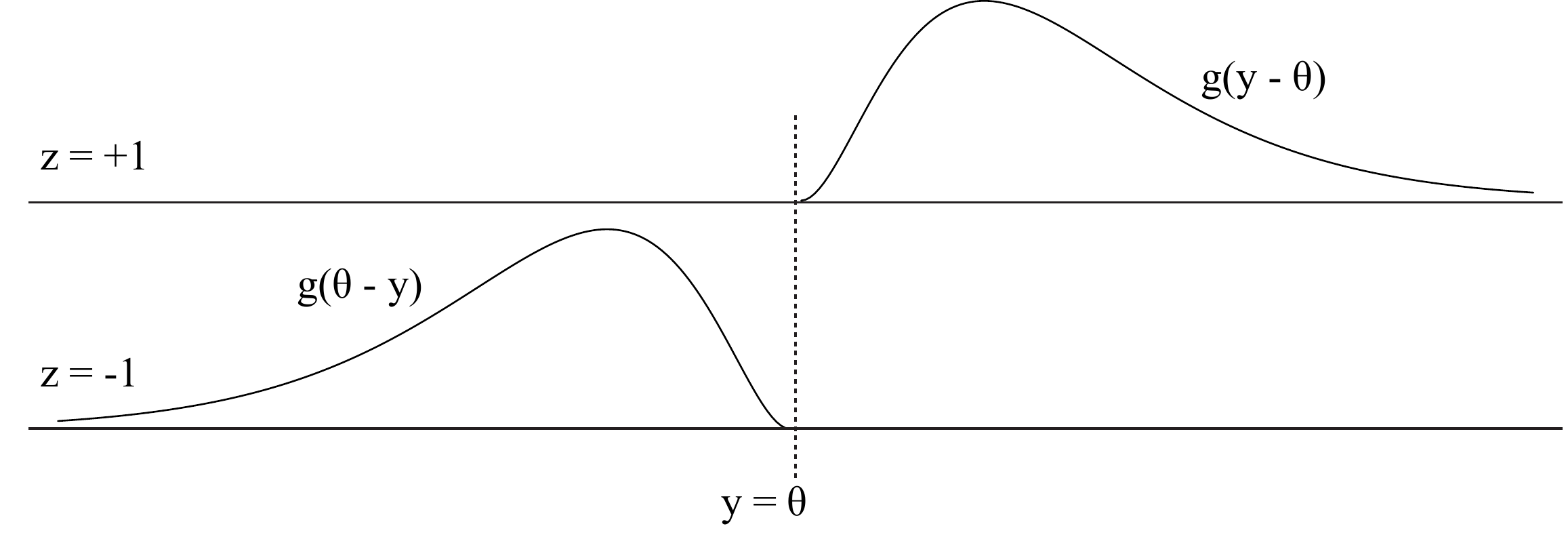}
\medskip

Here are the pertinent statistical facts for a single observation~$x=(y,z)$ 
from~$P_\th$. (See Section~\ref{DQM} for some proofs.) Define~$S(x)=y$
and~$A(x)=z$.
The marginal distribution $Q_\th$ of~$S$ has density
$$
h_\th(y) = \tfrac12g(y-\th)+\tfrac12 g(\th-y)
\qt{\wrt/ $\Leb$.}
$$

\begin{enumerate}
\item[(a)] 
The statistic~$A$
is ancillary (its distribution,~$\nu$, does not depend on~$\th$). By itself it gives no statistical information about~$\th$, but in conjunction with~$S$ it does tell us something about~$\th$: if $z=+1$ then (with $P_\th$ probability one) $\th<y$, and if $z=-1$ then~$y<\th$.

\item[(b)]
The statistic~$S$ is not sufficient because $\{z=+1\}=\{y>\th\}$ a.s.$[P_\th]$, implying $P_\th(Z=+1\mid y) = \{y>\th\}$ almost surely.  Equivalently,
$$
P_\th(A(x)=1\mid S)=\{S(x)>\th\}  
\qt{a.s.$[P_\th]$},
$$
which depends on~$\th$. (More formally, if $S$ were sufficient there would exist some measurable function~$\pi(y)$ for which $P_\th(z=1\mid S)=\pi(S(x))$ a.s.$[P_\th]$, for every~$\th$.)

\item[(c)]
Both~$\pp:=\{P_\th:\th\in\RR\}$ and~$\qq:=\{Q_\th:\th\in\RR\}$ have finite Fisher information: $\II_\qq(\th) = \II_\pp(\th) = \II$ for all $\th$,
with~$\II$ as in~(\ref{II.def}).
\end{enumerate}
\noindent
In short: There is no loss of Fisher information when only~$S(x)$ is observed, even though~$S$ is not a sufficient statistic.

\Remark 
K\&S used a slightly more involved construction, with density
\begin{align*}
f(x,\th) &= \{z=+1\}\left[0.7 g(y-\th) +0.3 g(\th-y)\right]
\\
&\quad+\{z= - 1\}\left[0.3 g(y-\th) +0.7 g(\th-y)\right]
\qt{where $x = (y,z)\in\xx$}
\end{align*}
\wrt/ $\Leb\otimes \mu$ where $\mu\{+1\}=\a = 1-\mu\{-1\}$ and $\a\ne 1/2$. The analysis that I present can be extended to this~$f_\th$.
\endRemark

\section{DQM interpretation}\label{DQM}
K\&S attributed the phenomenon in their version of the example in
Section~\ref{KSeg} to a failure of strict convexity of Fisher
information \wrt/ mixtures of statistical experiments. There is
another explanation involving the geometry of Hellinger derivatives,
which I find more illuminating.

By a theorem of \citet[Lemma A.3]{Hajek:72berk}, Lebesgue
integrability of the function~$\gdot^2/g$ in~\cref{II.def} implies that the  set
of densities $\gg:=\{g(y-\th):\th\in \RR\}$ (\wrt/ Lebesgue measure) is Hellinger
differentiable with Hellinger derivative $\gam(y-\th)$ at~$\th$, where
$$
\gam(w) := \frac{-\gdot(w)}{2\sqrt{g(w)}} = \frac{(2-w)}{2\sqrt{2}}e^{-w/2}\{w>0\}  .
$$
That is,
$$
\int\left|\sqrt{g(y-\th-t)} -\sqrt{g(y-\th)} - t\gam(y-\th)\right|^2dy = o(t^2)
\qt{as $t\to 0$}.
$$
This assertion is also easy to check by explicit calculations.
(See \citet[Cor.~12.2.1]{LehmannRomano2005} for details.)

The family of densities $\ff:=\{f_\th(x): \th \in\RR\}$, for $f_\th$ as in~\cref{f.def}, inherits
the Hellinger differentiability from~$\gg$:
\begin{equation}
\int\left|\sqrt{f_{\th+t}(x)} -\sqrt{f_{\th}(x)} - t\zeta_\th(x)\right|^2
\lam(dx)  = o(t^2)
\qt{as $t\to 0$},
\label{f.deriv}
\end{equation}
for the Hellinger derviative
$$
\zeta_\th(x) := \{z= +1\}\gam(y-\th) -\{z= -1\}\gam(\th-y).
$$

The significance of approximation~\cref{f.deriv}  becomes clearer when it is rewritten as a differentiability property of the likelihood ratios.
That is, it helps to work with the square root of the density
of~$P_{\th+t}$ \wrt/~$P_\th$. Unfortunately, $P_{\th+t}$ is not
dominated by~$P_\th$. In general, to eliminate such an embarrassment
one needs to split~$P_{\th+t}$ into a singular part~$\singP_{t,\th}$,
which concentrates on a set of zero~$P_\th$ measure, plus a part~$P_{\th+t}^{(abs)}$ that
has a density~$p_{t,\th}$ \wrt/~$P_\th$. For reasons related to the
asymptotic theory for repeated sampling, it is customary to make a
small extra assumption about the behavior of~$\singP_{t,\th}\xx$
as~$t$ tends to zero. Following \citet[Section~17.3]{LeCam:86book} and
\citet[Section~7.2]{LeCam:Yang:2000book}, I will call the slightly
stronger property \newdef{differentiability in  quadratic mean (DQM)},
to stress that the definition requires a little more than Hellinger
differentiability. 

\Remark 
Some authors (for example,
\citealt[page~457]{BKRW:93book}) use the term DQM as a synonym for
Hellinger differentiability.
\endRemark 

\begin{definition}	\label{dqm.def}
Say that $\pp=\{P_\th:\th \in \Th\}$, with $\Th\subseteq\RR$, is differentiable in quadratic mean (DQM) at~$\th$ with score function~$\Del_\th(x)$
  if, for $\th+t\in\Th$,
\begin{enumerate} 
\item
for the part~$\singP_{t,\th}$ of~$P_{\th+t}$ that is singular \wrt/~$P_\th$,
$$
\singP_{t,\th}(\xx)=o(t^2) 
\qt{as $t\to0$}
$$

\item
$\Del_\th\in\ll^2(P_\th)$  

\item
the absolutely continuous part of $P_{\th+t}$ 
has density $p_{t,\th}(x)$ \wrt/ $P_\th$ for which
$$
\sqrt{p_{t,\th}(x)} = 1 + \tfrac12 t \Del_\th(x) + r_{t,\th}(x)
\qt{with $P_\th \left(r_{t,\th}^2\right) =
o(t^2)$ as $t\to0$}.
$$
\end{enumerate} 
Call $\pp$  DQM if it is DQM at each~$\th$ in~$\Th$.
\end{definition}

\Remark 
The factor of~$1/2$ in requirement~(iii) ensures that $P_\th \Del_\th^2$ is equal to the Fisher information~$\II_\pp(\th)$ if the densities are suitably smooth in a pointwise sense. 
\endRemark

The~$\pp$ from Section~\ref{KSeg} is, in fact, DQM. For $t>0$ the
singular part~$\singP_{t,\th}$ has density
$\{z=-1\}g(\th-y)\{\th<y<\th+t\}$ \wrt/~$\lam$, so that
$\singP_{t,\th}(\xx)=O(|t|^3)$. The part of~$P_{\th+t}$ that is
dominated by~$P_\th$ has density
\begin{align}
p_{t,\th}(x) &= \frac{f_{\th+t}(x)}{f_\th(x)}\{f_\th(x)>0\}
\notag
\\
&= \{z = +1\}\frac{g(y-\th-t)}{g(y-\th)} \{y>\th\}
+ \{z = -1\}\frac{g(\th+t-y)}{g(\th-y)} \{y <\th\} .
\label{ptth}
\end{align}
There is a similar expression for the case $t<0$.
The score function equals
\begin{align}
\Del_\th(x) 
&= 2\frac{\zeta_\th(x)}{\sqrt{f_\th(x)}} \{f_\th(x)>0\}
\notag\\
&= \{z= +1\}\frac{\gam(y-\th)}{\sqrt{g(y-\th)}}\{y>\th\}
- \{z= -1\}\frac{\gam(\th-y)}{\sqrt{g(\th-y)}}\{y<\th\} .
\label{scoreP}
\end{align}

The density~$p_{t,\th}$ and the score function $\Del_\th(x)$ are
uniquely determined only up to a~$P_\th$ equivalence. As noted for fact~(b)
near the end of Section~\ref{KSeg},  sufficiency fails for~$S$ because
\begin{equation*}
\{z=+1\}=\{y>\th\}
\qt{a.s.$[P_\th]$}
\end{equation*}
Similarly $\{z=-1\} =\{y<\th\}$ a.s.$[P_\th]$. Together these two equivalences  explain why no Fisher information is lost. The score function~$\Del_\th$ is
changed only on a $P_\th$-negligible set if we omit the two indicator
functions involving~$z$ from~\cref{scoreP}. In effect, the score
function~$\Del_\th(x)$ depends on~$x$ only through the value of the
statistic~$S$. As the next theorem (which is proved in Section~\ref{Proof})
shows, that property is equivalent to the preservation of  Fisher information.

\Theorem  \label{preserveDQM}
Suppose $\pp=\{P_\th:\th\in\Theta\}$ on~$(\xx,\aa)$
is DQM  with score function~$\Del_\th$.
Suppose $S$ is a measurable map from $(\xx,\aa)$ into $(\yy,\bb)$
and $Q_\th = SP_\th$ is the distribution of~$S$ under~$P_\th$.
Then: 
\begin{enumerate} 
\item
The statistical experiment
$\qq=\{Q_\th:\th\in\Th\}$ is also DQM, with score function 
$\tDel_\th(y) = P_\th(\Del_\th\mid S=y)$.  

\item
At each fixed~$\th$, Fisher information is preserved (that is, $\II_\pp(\th)=\II_\qq(\th)$) if and only if $\Del_\th(x) = \tDel_\th(Sx)$ \almev{P_\th}.

\end{enumerate} 
\endTheorem

With only notational changes, the  Theorem extends to the case
where~$\Th$ is a subset of some Euclidean space; no extra conceptual
difficulties arise in higher dimensions.

\subheading{Credit where credit is due.}
The results stated in Theorem~\ref{preserveDQM} have an interesting history. Property~(i) was asserted (``Direct calculations show that the function~$q^{1/2}(y;\th)$ is differentiable in $L_2(\widetilde\nu)$ and possess a continuous derivative \dots") in Theorem~7.2 of 
\citet[Chapter~I, page~70]{Ibragimov:Hasminskii:81book}, an English translation from the 1979 Russian edition. However, that Theorem also (incorrectly, as noted by K\&S) asserted that Fisher information is preserved if and only if~$S$ is sufficient.

\citet[pages~19--21]{Pitman1979book} established differentiability in mean, a property slightly different from~(i), in order to deduce a result equivalent to~(ii).

\citet[Section~7]{Lecam:Yang:88AS} deduced an  analogue of~(i)  (preservation
of DQM under restriction to sub-sigma-fields) by an indirect argument using
equivalence of DQM with the existence of a quadratic approximation to
likelihood ratios of product measures (an LAN condition). 

\citet[page~461]{BKRW:93book} proved result~(i), citing 
\citet{Ibragimov:Hasminskii:81book},
\citet{Lecam:Yang:88AS}, and 
\citet[Appendix~A3]{vaart:88cwi} for earlier proofs. The last of these was a revised
(``I have not resisted the temptation to rewrite numerous parts of the original manuscript") version of \citeauthor{vaart:88cwi}'s 1987 Ph.D. thesis. He  cited \citet{Lecam:Yang:88AS} and a manuscript version of \citet{BKRW:93book}.

\section{Large sample interpretation}
The example in Section~\ref{KSeg} shows that, for a sample~$x=(y,z)$ of size one from~$P_\th$, some ``statistical information about~$\th$ (namely, whether~$\th<y$ or~$\th>y$) is lost if we discard~$z$. The loss is not detected by Fisher's measure of information. An asymptotic analysis, based on a sample of size~$n$ from~$P_\th$, sheds a little light on why the~$z$ contribution is relatively unimportant.

Write $\PP_{\th,n}$ and~$\QQ_{\th,n}$ for the $n$-fold product measures~$P_\th^n$
and~$Q_\th^n$, with the probability measures~$P_\th$ and $Q_\th$ as in Section~\ref{KSeg}.  That is, the statistical experiment
$\pp_n=\{\PP_{\th,n}: \th\in \Th\}$ corresponds to taking~$n$ independent observations~$x_1=(y_1,z_1),\dots, x_n=(y_n,z_n)$ from~$P_\th$ and~$\qq_n=\{\QQ_{\th,n}:\th\in\Th\}$ corresponds to~$y_1,\dots,y_n$. 

Both~$\pp_n$ and~$\qq_n$ are locally asymptotically normal
\cite[Chapter~6]{LeCam:Yang:2000book}. They share the same local normal approximations because that have the same score functions and (hence) the same Fisher information functions: for each fixed~$\th$ and each finite subset~$T$ of the real line, the ``local experiments''
\begin{equation*}
\{\PP_{\th+tn^{-1/2},n}:t\in T\}\qt{and}\qquad
\{\QQ_{\th+tn^{-1/2},n}:t\in T\}
\end{equation*} 
are asymptotically equivalent in the  (weak) Le~Cam sense. The deficiency distance
\cite[Section~2.2]{LeCam:Yang:2000book} between these two local experiments tends to zero as~$n$ tends to infinity. 
The local asymptotic equivalence of~$\pp_n$ and~$\qq_n$ has many  consequences. For example, classical theory establishes existence of many different estimators~$\that=\that(x_1,\dots,x_n)$ for which $\sqrt{n}\left(\that-\th\right)$ converges in distribution under~$\PP_{\th,n}$ to~$N(0,\II^{-1})$, and  
many different estimators~$\th_n^*= \th_n^*(y_1,\dots,y_n)$ for which $\sqrt{n}\left(\th_n^*-\th\right)$ converges in distribution under~$\QQ_{\th,n}$ to the same~$N(0,\II^{-1})$.
As shown by the H\'ajek-Le~Cam convolution and asymptotic minimax theorems \cite[Section~2.3]{BKRW:93book}, there are various senses in which 
the~$N(0,\II^{-1})$ limit is the best we can hope to achieve for 
either experiment.
Asymptotically speaking, the $z_i$'s must be contributing at a less important level. 

\Remark
Except for the purpose of the root-$n$ asymptotics,  perhaps we should agree with \citet[Section~5]{Basu1975StatInfo} that  the Fisher information function is a ``mathematically interesting but statistically rather fruitless notion".
\endRemark

For $i=1,\dots,n$, write~$\yL$ for the largest $y_i$   for which $z_i = -1$  and $\yR$ the smallest~$y_i$ for which~$y_i=+1$. Each~$z_i$ tells us whether~$y_i>\th$ or~$y_i<\th$, implying
\begin{equation}
\yL<\th <\yR\qt{with $\PP_{\th,n}$ probability one}.
\label{z.int}
\end{equation}

The~$w^2$ decay in~$g(w)$ at zero, implies that both $\th-\yL$ and~$\yR-\th$ are decreasing at an~$n^{-1/3}$ rate. In fact both~$n^{1/3}(\th-\yL)$ and $n^{1/3}(\yR-\th)$ have nontrivial limit distributions under~$\PP_{\th,n}$. For example, for each~$s>0$ direct calculation shows that 
$P_\th(\th,\th+sn^{-1/3})=s^3/(6n)+ o(1/n)$, so that
\begin{equation*}
\PP_{\th,n}\{n^{1/3}(\yR-\th)>s\}
=\PP_{\th,n}\{\text{no $z_i$'s in $(\th,\th+sn^{-1/2})$}\} \to \exp(-s^3/6).
\end{equation*}
For any $n^{-1/2}$-consistent estimator~$\that$
the event~$A_n=\{\yL<\that<\yR\}$ has $\PP_{\th,n}$ probability that tends very rapidly to one. Put another way,
\begin{equation*}
\PP_{\th,n}\{\exists i\le n: \th < y_i <\that\text{ or }\that<y_i <\th\}\to0 .
\end{equation*}
With high probability, what we learn from the~$z_i$'s just duplicates what we usually can learn from the~$y_i$'s.

To make the idea more concrete,
define $\zin = \text{sgn}(y_i-\that)$ and $\xin = (y_i,\zin)$.  That is, 
\begin{equation*}
\zin =
\Cases
+1&if $y_i >\that$
\\
-1& if $y_i < \that$
\endCases
.
\end{equation*}
On the event~$A_n$ we have $x_i=\xin$ for $i=1,\dots,n$. If $\PP^*_{\th,n}$ denotes the joint distribution of $x^*_{1,n},\dots,x^*_{n,n}$ then
$$
\supnl_{\th\in\Th}\normTV{\PP^*_{\th,n}-\PP_{\th,n}}
\le \supnl_{\th\in\Th}\PP_{\th,n} A_n^c\to 0
.
$$
In the terminology of Le~Cam's theory for convergence of statistical experiments,~$\pp_n$ and~$\qq_n$ are asymptotically equivalent, not just locally asymptotically equivalent in the weak sense. 
The vector $(y_1,\dots,y_n)$ is asymptotically sufficient for~$\pp_n$, in Le~Cam's sense. The map $(y_1,\dots,y_n)\mapsto(x^*_{1,n},\dots,x^*_{n,n})$
defines a Le~Cam transition \cite[Theorem~2.2]{LeCam:Yang:2000book} that bounds the deficiency~$\del(\qq_n,\pp_n)$.

Put another way, for every statistic~$\psi_n(x_1,\dots,x_n)$ for~$\pp_n$ there is another statistic $\psi_n^*(y_1,\dots, y_n) =\psi_n(x^*_{1,n},\dots,x^*_{n,n})$ for $\qq_n$ that has the same asymptotic behavior.

\Remark 
Rough calculations suggest that the Le~Cam distance between~$\pp_n$ and~$\qq_n$ tends to zero like $\exp(-Cn^{1/3})$ for some constant~$C$. I omit the details because the actual rate is not important for the story I am telling.
\endRemark

\section{Proof of Theorem~\protect{\ref{preserveDQM}}}  \label{Proof}
Recall  that the Kolmogorov conditional expectation
$P_\th(\cdot\mid S=y)$ is abstractly defined, via the
Radon-Nikodym theorem, as an increasing linear map (depending on~$\th$)
$\kappa:L^1(P_\th)\to L^1(Q_\th)$ with properties analogous to those
enjoyed by a Markov kernel. If we identify an~$f$
in~$L^1(P_\th)$ with the (signed) measure~$\mu_f$ for which
$d\mu_f/dP_\th=f$, then $g=\kappa f$ is the density of~$S\mu_f$
\wrt/~$Q_\th$. 
To stress the analogy with Markov kernels I will write $\kappa_y f$, or
even $\kappa_y f(x)$, instead of~$(\kappa f)(y)$. 
Thus the defining property of~$\kappa$ can be rewritten as
\begin{equation}
Q_\th f_1(y)\kappa_y f_2 = P_\th f_1(Sx)f_2(x)
\label{Kol}
\end{equation}
for measurable real functions $f_1$ on~$\yy$ and $f_2$ on~$\xx$,
at least when $f_1(Sx)f_2(x)$ is $P_\th$-integrable.
A reader
who chose to interpret~$\kappa_y$ as a Markov kernel would lose only a
tiny amount of generality.

Of course if one regards $\kappa$ as acting on the function~$\ll^1(P_\th)$, instead
of on the space~$L^1(P_\th)$ of $P_\th$-equivalence classes, then one
should qualify assertions with the occasional \almev{P_\th} caveats
and regard~$\kappa f$ as being defined only up to $Q_\th$ equivalence.
Following the usual custom, I will omit such qualifiers.

\subheading{Proof of assertion~(i).}
The following argument is adapted from
\citet[Appendix~A3]{vaart:88cwi}.

To simplify notation, I will prove that~$\qq$ is DQM only at~$\th=0$,
writing $\singP_t$ instead of~$\singP_{t,0}$ and $p_t$ instead
of~$p_{t,0}$. Keep in mind that~$\kappa_y$ now denotes the conditional expectation
operator~$P_0(\cdot\mid S=y)$. 
For each function~$h(x)$ in~$\ll^2(P_0)$ I will write $\htil(y)$ for its conditional mean~$\kappa_yh(x)$ and
$$
\var_y h := \kappa_y\left(h(x)-\htil(y)\right)^2 = \kappa_y h(x)^2 - \htil(y)^2
$$
for its conditional variance.

Start with the simplest case where~$P_t$ is actually dominated by~$P_0$. Then
$$
\xi_t(x) = \sqrt{{dP_t}/{dP_0}} = 1 +\tfrac12 t\Del_0(x) + r_t(x)
\qt{with $P_0r_t^2 = o(t^2)$}
$$
and
\begin{equation}
\xitil_t(y) := \kappa_y\xi_t(x) = 1+ \tfrac12 t \Deltil_0(y) + \rtil_t(y)
\qt{with $Q_0\rtil_t^2 \le P_0 r_t^2= o(t^2)$}.
\label{xitil}
\end{equation}
and, by the Radon-Nikodym property,
$$
\eta_t(y) =\sqrt{{dQ_t}/{dQ_0}} = \sqrt{\kappa_y\xi_t(x)^2}
\quad.
$$

The proof of assertion~(i) will work by showing that the difference $\del_t(y) := \eta_t(y)-\xitil_t(y)$ is small, in the sense that~$Q_0\del_t^2=o(t^2)$. For then we will have
$$
\eta_t(y) = 
1 +\frac12 t \Deltil_0(y) +\bigl[\rtil_t(y) +\del_t(y)\bigr]
\qt{with }Q_0\bigl[\rtil_t(y) +\del_t(y)\bigr]^2=o(t^2),
$$
which implies DQM for~$\qq$ at~$0$.

The desired property for~$\del_t$ will be derived from the following three facts about the conditional variance
\begin{equation}
\sig_t^2(y) := \var_y(\xi_t) = \kappa_y \xi_t(x)^2 - \xitil_t(y)^2= \eta_t(y)^2 - \xitil_t(y)^2.
\label{sigt}
\end{equation}
\begin{enumerate} 
\item[(a)]
The  representation $\sig_t^2(y) =\kappa_y\left(\xi_t(x)-\xitil_t(y)\right)^2$ gives
\begin{align*}
\sig_t^2(y)
&= \kappa_y\left(\tfrac12 t\left[\Del_0(x)-\Deltil_0(y)\right] +\left[r_t(x)-\rtil_t(y)\right]\right)^2
\notag\\
&\le 2\left(\tfrac12 t\right)^2\kappa_y\left[\Del_0(x)-\Deltil_0(y)\right] ^2 + 2\kappa_y\left[r_t(x)-\rtil_t(y)\right]^2
\notag\\
&\le \tfrac12 t^2 \kappa_y\Del_0^2 + 2\kappa_y r_t^2   .
\end{align*}
\end{enumerate} 

\Remark 
The cancellation of the leading $1$ when $\xitil_t$ is subtracted from~$\xi_t$ seems to
be vital to the proof. For general Hellinger differentiability, the cancellation would not occur.
\endRemark

\begin{enumerate} 

\item[(b)]
 $\del_t(y) \ge 0$ because
 $\eta_t(y)^2-\xitil_t(y)^2 =\sig_t^2(y)\ge0$.
\smallskip
\item[(c)]
Substitution of $\del_t+\xitil_t$ for~$\eta_t$ in~\cref{sigt} gives
$$
\sig_t^2(y) = 2\del_t(y)\xitil_t(y) + \del_t(y)^2 .
$$
\end{enumerate} 

The rest is easy.
For each $\eps>0$ define
$$
\Ate :=\{y\in \yy: \xitil_t(y) \ge\tfrac12,\, \sig_t(y) \le \eps\}.
$$
Integration of inequality~(a) gives
$$
Q_0 \sig_t^2(y) \le \tfrac12 t^2 P_0\Del_0^2 + 2P_0r_t^2 = O(t^2)+o(t^2)\le Ct^2
\qt{for some constant $C$},
$$
which, together with~\cref{xitil}, implies $Q_0\Ate\to1$ as~$t\to0$.

On the set~$\Ate$ equality~(c) ensures that $\del_t(y) \le \sig_t^2(y) \le \eps \sig_t(y)$; on $\Ate^c$ the nonnegativity of~$\del_t$ and equality~(c) give~$\del_t^2\le \sig_t^2$.
Thus
\begin{align*}
Q_0 \del_t(y)^2
&\le \eps^2 Q_0 \sig_t^2(y)\{y\in \Ate\} 
+ Q_0 \sig_t^2(y)\{y\notin \Ate\}
\\
&\le \eps^2 Ct^2 + \tfrac12 t^2 Q_0 \kappa_y\Del_0^2\Ate^c + 2Q_0\kappa_y r_t^2
\qt{by (a)}.
\end{align*}
The $Q_0$-integrability of $\kappa_y\Del_0^2$ and the Dominated Convergence theorem imply $Q_0 \kappa_y\Del_0^2\Ate^c\to0$. It follows that~$Q_0\del_t^2=o(t^2)$.

Finally, what happens when~$P_t$ is not dominated by~$P_0$?
The analysis for $\xi_t^2$, the density of the part of~$P_t$ that is dominated by~$P_0$, is the same as before. The image measure~$S\singP_t$ has total mass of order~$o(t^2)$, part of which gets absorbed into~$\singQ_t$. The part of~$S\singP_t$ that is dominated by~$Q_0$ contributes an extra nonnegative term,~$\gam_t(y)$, to the density~$dQ_t^{(abs)}/dQ_0$.  The $\eta_t^2(y)$ becomes $\kappa_y \xi_t^2(y) +\gam_t(y)$. The extra term causes little  trouble because
$$
\sqrt{\kappa_y\xi_t^2} 
\le \eta_t \le
\sqrt{\kappa_y\xi_t^2} +\sqrt{\gam_t}
\qt{and}\qquad Q_0\gam_t=o(t^2).
$$

\subheading{Proof of assertion~(ii).}
Write $\HH$ for the closed subspace of~$L^2(P_\th)$ consisting of (equivalence classes of) functions measurable \wrt/ the sub-sigma-field of~$\aa$ generated by~$S$. Each member of~$\HH$ is of the form $f(Sx)$ for some~$f$ in~$L^2(Q_\th)$. The orthogonal projection of~$\Del_\th$ onto~$\HH$ equals~$\Deltil_\th(Sx)$. Thus
$$
\II_\pp(\th) = P_\th \Del_\th(x)^2 = 
P_\th \Deltil_\th(Sx)^2 + P_\th\left[\Del_\th(x)-\Deltil_\th(Sx)\right]^2.
$$
The first term on the right-hand side equals $Q_\th(\Deltil_\th^2)= \II_\qq$; the last term is zero if and only if $\Del_\th(x)=\Deltil_\th(Sx)$ \almev{P_\th}.

\section*{Acknowledgements}
Many thanks to the referee for pointing out the connections with the work of Basu.
In retrospect, it is a little surprising that I was unable to find a counterexample analogous to that of K\&S in the volume \citep{DasGupta2011} of Basu's selected works.

\bibliographystyle{imsart-nameyear}%
\bibliography{DBP}%

\end{document}